\documentclass[a4paper,11pt]{article}
\usepackage{latexsym}
\usepackage{amssymb}
\usepackage{amsfonts}
\usepackage{amsmath}
\usepackage{mathabx}
\usepackage{indentfirst}
\usepackage{graphicx}
\usepackage{epsfig}

\usepackage{fullpage}

\usepackage{color}
\usepackage{mathrsfs}
\usepackage{manfnt}

\usepackage{tikz,ifthen}

\newtheorem{prop}{Proposition}[section]
\newtheorem{teor}{Theorem}[section]
\newtheorem{lemma}{Lemma}[section]
\newtheorem{cor}{Corollary}[section]

\newcommand{\vDashv}{%
  \mathrel{%
    \text{%
      \ooalign{$\vDash$\cr$\Dashv$\cr}%
    }%
  }%
}
\newcommand{\cvd}{\hfill $\blacksquare$\bigskip}
\date{}

\author{Luca Ferrari\thanks{Dipartimento di Matematica e Informatica ``U. Dini", viale Morgagni 65, 50134 Firenze, Italy
{\tt luca.ferrari@unifi.it}.\newline \indent Partially supported
by MIUR PRIN 2010-2011 grant ``Automi e Linguaggi Formali: Aspetti
Matematici e Applicativi", code H41J12000190001, and by INdAM-GNCS 2014 project ``Studio
di pattern in strutture combinatorie".}}

\title{Dyck algebras, interval temporal logic and posets of intervals}

\begin{document}

\maketitle

\begin{abstract}

We investigate a natural Heyting algebra structure on the set of Dyck paths of the same length.
We provide a geometrical description of the operations of pseudocomplement and relative pseudocomplement,
as well as of regular elements. We also find a logic-theoretic interpretation of such Heyting algebras,
which we call \emph{Dyck algebras}, by showing that they are the algebraic counterpart
of a certain fragment of a classical interval temporal logic (also known as Halpern-Shoham logic).
Finally, we propose a generalization of our approach, suggesting a similar study of the Heyting algebra
arising from the poset of intervals of a finite poset using Birkhoff duality. In order to illustrate this,
we show how several combinatorial parameters of Dyck paths can be expressed in terms of the
Heyting algebra structure of Dyck algebras together with a certain total order on the set of atoms
of each Dyck algebra.

\end{abstract}

\section{Introduction}

Among the plethora of different logics generalizing and extending the classical one,
a family of logics which has proved very useful especially in computer science is that of \emph{temporal logics}.
A temporal logic is essentially a kind of logic which allows to deal with statements whose truth values can vary in time.
Applications in computer science concern, for example, formal verification,
where temporal logics show their expressiveness in stating requirements of hardware or software systems.
Starting from the generic idea stated above, one can conceive several different types of temporal logics,
depending on the structure of time states and on how time states are managed.
A particularly interesting class of temporal logics are the so-called \emph{interval temporal logics}.
An interval temporal logic is characterized by the fact that the truth of a statement depends on
the \emph{time interval} it is evaluated on (rather than the time instant). Such kinds of logics are useful
when it is important to work with properties which remain true (or false) for a certain amount of time.
The relevance of these logics for computer science is even more evident: think, for instance, of processes,
for which it is meaningful to reason in terms of time intervals rather than time instants.
More generally, interval temporal logics have been successfully applied to temporal databases,
specification, design and verification of hardware components and concurrent real-time processes;
see, for instance \cite{GMS} and the references therein.

To work with any interval temporal logic, it is important to understand
which kinds of relations among intervals of time instants are relevant to the specific logic one wish to consider.
The classification of all possible such relations has been pursued by Allen \cite{AF},
who has also defined an algebraic structure to deal with them.
The modal logic of time intervals resulting by considering the whole of Allen's relations is usually referred to as
the \emph{Halpern-Shoham logic} \cite{HS}. Typically, one selects a subset of Allen's relations,
thus defining the related fragment of the Halpern-Shoham logic. Most studied in this context are decidability questions,
as witnessed by many works appeared in recent years (an example related to a fragment which is relevant to our paper
is \cite{MM}).

In the present paper we propose a combinatorial description of a specific interval temporal logic
whose underlying time model is a finite linear order. Specifically, we consider what is sometimes called
the \emph{logic of subintervals}, that is the interval temporal logic in which,
from the truth of a statement on a certain interval of time instants $I$,
the truth of that statement on all subintervals of $I$ follows.
We show that, given a linearly ordered set of time instants of cardinality $n-1$,
the algebraic counterpart of the associated logic of subintervals is given by
a certain Heyting algebra structure on the set of Dyck paths of semilength $n$,
which is more precisely the canonical Heyting algebra structure associated with the distributive lattice structure
on Dyck paths of semilength $n$ induced by ordering them by \emph{geometric inclusion}
(i.e. a Dyck path $P$ is declared to be less than or equal to a Dyck path $Q$ whenever,
in the usual two-dimensional drawing of Dyck paths, $P$ lies weakly below $Q$, see \cite{FP} and also the next section).
We also give a fully geometric description of relative pseudocomplement and pseudocomplement in such Dyck algebras,
thus supplementing similar results that have been illustrated in a more algebraic fashion in \cite{Muh}.
Finally, we try to broaden the scope of our work, by proposing a possible generalization.
The idea is to consider the poset of intervals (ordered by inclusion) of any poset $\mathcal{P}$
(rather than a totally ordered set) and to investigate properties of the Heyting algebra $\mathcal{H}$
obtained from $\mathcal{P}$ by classical (generalized) Birkhoff duality. More specifically, we ask what properties of $\mathcal{H}$ can be expressed
in terms of the partial order $\mathcal{P}$. In the specific case of a finite totally ordered $\mathcal{P}$
(which is the case studied in the present paper), we illustrate the above project from a combinatorial point of view,
namely we express several statistics of combinatorial interest in terms of the Heyting algebra structure of Dyck paths
together with the partial order structure on the atoms of such an algebra.
We close our paper by proposing some further directions of future research.

\section{Heyting algebras of Dyck paths}\label{structure}

%
%
%
%
%
%


Given a Cartesian coordinate system, a \emph{Dyck path} is a
lattice path starting from the origin, ending on the $x$-axis,
never falling below the $x$-axis and using only two kinds of
steps, $u(p)=(1,1)$ and $d(own)=(1,-1)$. A Dyck path can be
encoded by a word $w$ on the alphabet $\{ u,d\}$ such that in
every prefix of $w$ the number of $u$ is greater than or equal to
the number of $d$ and the total number of $u$ and $d$ in $w$ is
the same (the resulting language is called \emph{Dyck language}
and its words \emph{Dyck words}). The \emph{length} of a Dyck path
is the length of the associated Dyck word (which is necessarily an
even number). A \emph{peak} in a Dyck path is a pair of consecutive
steps of the form $ud$; a \emph{hill} is a peak at height 0
(i.e. lying on the $x$-axis). A \emph{factor} of a Dyck path is
any minimal subsequence of consecutive steps starting and ending on the
$x$-axis; every Dyck path can be clearly decomposed in a unique way
as the product (juxtaposition) of its factors. In particular,
a hill is also called a \emph{trivial factor}. A \emph{pyramid} is
a subsequence of consecutive steps of the form $u^k d^k$ ($k\geq 1$)
starting and ending on the $x$-axis. In particular, a hill is a pyramid.
A \emph{return} is a point of the path,
other than the starting one, lying on the $x$-axis.
We will usually refer to a return by using its abscissa
(which is necessarily an even number).

The set $D_n$ of Dyck paths of semilength $n$ can be endowed with
a very natural poset structure. Given $P,Q\in D_n$, we say that
$P\leq Q$ when, in the above described two-dimensional drawing of
Dyck paths, $P$ lies weakly below $Q$. Properties of the posets
$\mathcal{D}_n =(D_n ,\leq )$ have been investigated in
\cite{FM1,FM2,FM3,FP}. In particular, it is shown that
$\mathcal{D}_n$ is a distributive lattice, and for this reason it
will be called the \emph{Dyck lattice of order} $n$. We point out
that this last assertion is a consequence of the (easy to
observe) fact that $\mathcal{D}_n$ is isomorphic to the dual of the
Young lattice of integer partitions whose Ferrers diagrams fit
into the staircase diagram $(n-1,n-2,\ldots ,2,1)$ \cite{S}. The
language of Dyck paths, however, gives a geometric flavor to the
subject which allows to express several properties in a more
fascinating way, as well as to suggest possible analogies with
other families of lattice paths.

Recall that a \emph{join-irreducible element} of a poset $\mathcal{P}$ is an element $a$ such that,
if $a=x\vee y$, then $a=x$ or $a=y$. In particular, if $\mathcal{P}$ has minimum $0$,
an \emph{atom} is an element covering $0$ (hence an atom is join-irreducible).
Moreover, a subset $I$ of $\mathcal{P}$ is a \emph{down-set} whenever, for every $x,y$ in $\mathcal{P}$,
if $y\in I$ and $x\leq y$, then $x\in I$. The well-known \emph{Birkhoff representation theorem}
(see, for instance, \cite{DP}) states that every finite distributive lattice is isomorphic to
the lattice of down-sets of the poset of its join-irreducibles. As a consequence,
every element of a finite distributive lattice is the join of the join-irreducibles below it.
Concerning Dyck lattices, a join-irreducible is a path all of whose factors are hills
except for a single pyramid having at least 4 steps (see \cite{FM1}). In particular,
an atom is a join-irreducible in which the unique nontrivial pyramid has \emph{exactly}
4 steps.

Since Dyck lattices are finite distributive lattices, they also
have a canonical Heyting algebra structure. Recall that a
\emph{Heyting algebra} is a lattice $\mathcal{H}$ with minimum $0$
and maximum $1$ such that the relative pseudocomplement of $x$
with respect to $y$ exists for all $x,y\in \mathcal{H}$. By
definition, the \emph{relative pseudocomplement of $x$ with
respect to $y$} is the element $x\rightsquigarrow y$ defined as
follows:
$$
x\rightsquigarrow y =\bigvee \{ z\in \mathcal{H}\; |\; x\wedge z\leq y\} .
$$

The Heyting algebra of Dyck paths of semilength $n$ will be
denoted $\mathfrak{D}_n$, and we will call it the \emph{Dyck
algebra of order} $n$.

In a Heyting algebra $\mathcal{H}$, two important notions are
those of pseudocomplement and of regular element. The
\emph{pseudocomplement} of $x$ is defined as $\sim \! x=x\rightsquigarrow
0$. It can be shown that $x\leq \sim \sim \! x$. The converse,
however, does not hold in general. An element $x$ of $\mathcal{H}$
is said to be \emph{regular} whenever $x=\sim \sim \! x$. The
subposet of regular elements of a Heyting algebra forms a Boolean
algebra.

The main aim of the present section is to give a combinatorial
description of relative pseudocomplement and pseudocomplement in Dyck algebras,
as well as to characterize the Boolean algebra of the regular
elements. We point out that similar results have been obtained in
\cite{Muh}. Our statements, however, have a more geometric flavor,
which would hopefully result in a more natural way of capturing
the above mentioned notions.

\bigskip

For any pair of Dyck paths $(P,Q)$ of semilength $n$, 
we define the \emph{crossing set} $C(P,Q)\subseteq [2n]\cup \{ 0\}=\{ 0,1,2,\ldots ,2n\}$
of $(P,Q)$ by declaring $x\in C(P,Q)$ whenever exactly one of the following conditions holds:

\begin{enumerate}

\item $x\in \{ 0,2n\}$;

\item $P$ and $Q$ have a common point having abscissa $x$; moreover
$P$ has an up step starting at that point and $Q$ has a
down step starting at that point;

\item $P$ and $Q$ have a common point having abscissa $x$; moreover
$P$ has a down step arriving at that point and $Q$ has an
up step arriving at that point.

\end{enumerate}

Roughly speaking, an element of the crossing set of $(P,Q)$ is
either the abscissa of the starting/ending point of the two paths
or the abscissa of a point in which the two paths crosses in a specific way.
More precisely, suppose that $C(P,Q)=\{ x_0, x_1 ,x_2 ,\ldots ,x_k \}$,
where the $x_i$'s are listed in increasing order (so that $x_0 =0$
and $x_k =2n$). If $i$ is even, then $P$ lies weakly below $Q$
between $x_i$ and $x_{i+1}$ (``weakly" meaning that $P$ and $Q$
may coincide in some point other than those of abscissas $x_i$ and
$x_{i+1}$); if $i$ is odd, then $P$ lies strictly above $Q$
between $x_i$ and $x_{i+1}$. Notice that $k$ is necessarily an odd
number (or, which is the same, the cardinality of $C(P,Q)$ is
even): indeed, both at the beginning and at the end $P$ lies
weakly below $Q$ (since both paths necessarily starts with an up
step and ends with a down step). Finally, observe that clearly
$C(P,Q)\neq C(Q,P)$ in general.

\begin{prop} Let $P,Q\in D_n$ and let $C(P,Q)=\{ x_0, x_1 ,x_2 ,\ldots ,x_k
\}$ be the crossing set of $(P,Q)$. Then $P\rightsquigarrow Q\in D_n$
is the Dyck path constructed as follows:

\begin{enumerate}

\item if $i$ is even, then the portion of $P\rightsquigarrow Q$ between
$x_i$ and $x_{i+1}$ is the unique subpath of the form $u^\alpha
d^\beta$ whose starting and ending points are the same as $P$ and $Q$,
for suitable nonnegative integers $\alpha$ and $\beta$;

\item if $i$ is odd, then $P\rightsquigarrow Q$ coincides with $Q$
between $x_i$ and $x_{i+1}$.

\end{enumerate}

\end{prop}

\emph{Proof.}\quad We observe that, if $i\neq 0$ is even, then
necessarily $P$ has an up step starting at abscissa $x_i$ and $Q$
has a down step starting at abscissa $x_i$, whereas, if $i\neq k$
is odd, then $P$ has a down step ending at abscissa $x_i$ and $Q$
has an up step ending at abscissa $x_i$. Thus, between $x_i$ and
$x_{i+1}$, if $i$ is even then $P$ lies weakly below $Q$,
otherwise (i.e. if $i$ is odd) $P$ lies strictly above $Q$ (this
last statement is true also in the cases $i=0,k$). As a
consequence, if $i$ is even, $P\rightsquigarrow Q$ can run as high
as possible between $x_i$ and $x_{i+1}$; this is achieved by
putting as many up steps as possible immediately after $x_i$,
followed by the correct number of down steps, which means that the
portion of $P\rightsquigarrow Q$ between $x_i$ and $x_{i+1}$ is of the
form $u^\alpha d^\beta$, as required. On the other hand, if $i$ is
odd, then $P\rightsquigarrow Q$ must coincide with $Q$ between $x_i$
and $x_{i+1}$, in order to have $(P\rightsquigarrow Q )\wedge P\leq Q$,
and this is clearly the maximum subpath between $x_i$ and
$x_{i+1}$ which satisfies such a condition.\cvd

The result of the above proposition can be restated less formally,
but maybe more expressively, as follows:
$P\rightsquigarrow Q$ is obtained from $Q$ by replacing those portions of path in which
$P$ lies weakly below $Q$ with the highest possible Dyck factors.

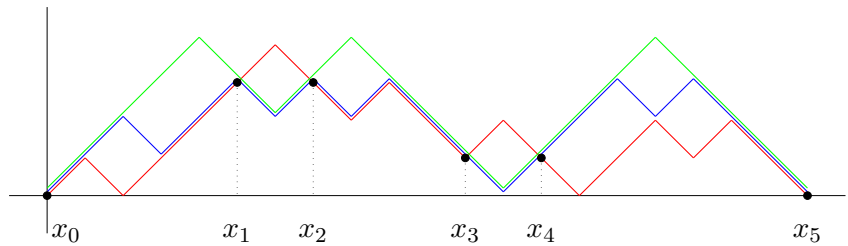
\begin{figure}[h!]
\begin{center}
\begin{tikzpicture}
\begin{scope}[scale=0.5]

\draw[very thin] (-1,0) -- (21,0);
\draw[very thin] (0,-1) -- (0,5);

\draw[red] (0,0) -- (1,1);
\draw[red] (1,1) -- (2,0);
\draw[red] (2,0) -- (6,4);
\draw[red] (6,4) -- (8,2);
\draw[red] (8,2) -- (9,3);
\draw[red] (9,3) -- (11,1);
\draw[red] (11,1) -- (12,2);
\draw[red] (12,2) -- (14,0);
\draw[red] (14,0) -- (16,2);
\draw[red] (16,2) -- (17,1);
\draw[red] (17,1) -- (18,2);
\draw[red] (18,2) -- (20,0);

\draw[blue] (0,0.1) -- (2,2.1);
\draw[blue] (2,2.1) -- (3,1.1);
\draw[blue] (3,1.1) -- (5,3.1);
\draw[blue] (5,3.1) -- (6,2.1);
\draw[blue] (6,2.1) -- (7,3.1);
\draw[blue] (7,3.1) -- (8,2.1);
\draw[blue] (8,2.1) -- (9,3.1);
\draw[blue] (9,3.1) -- (12,0.1);
\draw[blue] (12,0.1) -- (15,3.1);
\draw[blue] (15,3.1) -- (16,2.1);
\draw[blue] (16,2.1) -- (17,3.1);
\draw[blue] (17,3.1) -- (20,0.1);

\draw[green] (0,0.2) -- (4,4.2);
\draw[green] (4,4.2) -- (6,2.2);
\draw[green] (6,2.2) -- (8,4.2);
\draw[green] (8,4.2) -- (12,0.2);
\draw[green] (12,0.2) -- (16,4.2);
\draw[green] (16,4.2) -- (20,0.2);

\draw (0,0) [fill] circle (.1);
\draw (5,3) [fill] circle (.1);
\draw (7,3) [fill] circle (.1);
\draw (11,1) [fill] circle (.1);
\draw (13,1) [fill] circle (.1);
\draw (20,0) [fill] circle (.1);

\draw (0.5,-1) node {$x_0$};
\draw (5,-1) node {$x_1$};
\draw (7,-1) node {$x_2$};
\draw (11,-1) node {$x_3$};
\draw (13,-1) node {$x_4$};
\draw (20,-1) node {$x_5$};

\draw[dotted,very thin] (5,3) -- (5,0);
\draw[dotted,very thin] (7,3) -- (7,0);
\draw[dotted,very thin] (11,1) -- (11,0);
\draw[dotted,very thin] (13,1) -- (13,0);

\end{scope}
\end{tikzpicture}
\end{center}
\caption{$P$ is red, $Q$ is blue and $P\rightsquigarrow Q$ is green. \label{implication}}
\end{figure}

In Figure \ref{implication} we give an example of how to compute $P\rightsquigarrow
Q$ starting from $P$ and $Q$, as described in the above
proposition.

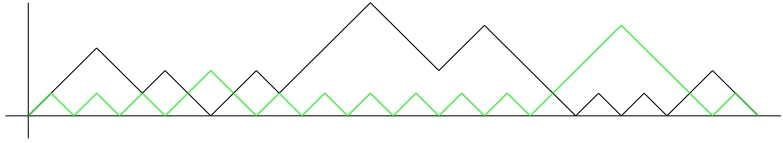
\begin{figure}[h!]
\begin{center}
\begin{tikzpicture}
\begin{scope}[scale=0.3]

\draw[very thin] (-1,0) -- (33,0);
\draw[very thin] (0,-1) -- (0,5);

\draw (0,0) -- (3,3);
\draw (3,3) -- (5,1);
\draw (5,1) -- (6,2);
\draw (6,2) -- (8,0);
\draw (8,0) -- (10,2);
\draw (10,2) -- (11,1);
\draw (11,1) -- (15,5);
\draw (15,5) -- (18,2);
\draw (18,2) -- (20,4);
\draw (20,4) -- (24,0);
\draw (24,0) -- (25,1);
\draw (25,1) -- (26,0);
\draw (26,0) -- (27,1);
\draw (27,1) -- (28,0);
\draw (28,0) -- (30,2);
\draw (30,2) -- (32,0);

\draw[green] (0,0) -- (1,1);
\draw[green] (1,1) -- (2,0);
\draw[green] (2,0) -- (3,1);
\draw[green] (3,1) -- (4,0);
\draw[green] (4,0) -- (5,1);
\draw[green] (5,1) -- (6,0);
\draw[green] (6,0) -- (8,2);
\draw[green] (8,2) -- (10,0);
\draw[green] (10,0) -- (11,1);
\draw[green] (11,1) -- (12,0);
\draw[green] (12,0) -- (13,1);
\draw[green] (13,1) -- (14,0);
\draw[green] (14,0) -- (15,1);
\draw[green] (15,1) -- (16,0);
\draw[green] (16,0) -- (17,1);
\draw[green] (17,1) -- (18,0);
\draw[green] (18,0) -- (19,1);
\draw[green] (19,1) -- (20,0);
\draw[green] (20,0) -- (21,1);
\draw[green] (21,1) -- (22,0);
\draw[green] (22,0) -- (26,4);
\draw[green] (26,4) -- (30,0);
\draw[green] (30,0) -- (31,1);
\draw[green] (31,1) -- (32,0);

\end{scope}
\end{tikzpicture}
\end{center}
\caption{A Dyck path (black) and its pseudocomplement (green). \label{pseudocomplement}}
\end{figure}

As a consequence, we have the following result, which gives us a recipe to compute
pseudocomplements in Dyck algebras (see Figure \ref{pseudocomplement}).
In the statement of the corollary, we will use the expression \emph{``sequence of $k$ consecutive hills"},
which should be clear in the case $k>0$. By convention, with the expression
\emph{``sequence of 0 consecutive hills"} we will mean a point of the path lying on the $x$-axis
(other than the starting and the ending ones) and neither preceded nor followed by a hill
(in other words, a return between two nontrivial factors).

\begin{cor} Let $P\in D_n$. Then $\sim \! P=P\rightsquigarrow 0$ is obtained from $P$ by
\begin{enumerate}
\item replacing each sequence of $k\geq 0$ consecutive hills starting at abscissa $x$ and ending at abscissa $x'$
with a pyramid of suitable height starting at $\max (0,x-2)$ and ending at $\min (x'+2,2n)$, and

\item completing the path by suitably adding a (finite) set of hills.

\end{enumerate}
\end{cor}

To conclude this section, we will give a characterization of regular elements of Dyck algebras.
Similarly to the previous results, our description will be in terms of the geometric shape of
the path.

\begin{prop} A Dyck path is regular if and only if its factors are all pyramids.
\end{prop}

\emph{Proof.}\quad For any Dyck path $P$, it follows from the previous corollary that
all factors of $\sim \! P$ are pyramids. Therefore, if $P$ is regular, then $P=\sim \sim \! P$,
and all factors of $P$ are pyramids.

For the converse, observe that the pseudocomplement operation
exchanges returns and non-returns of a Dyck path (that is, $(x,0)$ is
a return of $P$ if and only if $(x,0)$ is not a return of $\sim \! P$).
Now, if $P$ is a concatenation of pyramids, then $P$ is uniquely
determined by its returns, and the above observation implies that
$\sim \sim \! P=P$, i.e. $P$ is regular.\cvd

Recall that, given a poset $\mathcal{P}$, a \emph{closure operator} is a map
$\overline{\phantom{a}}:\mathcal{P}\rightarrow \mathcal{P}$ such that, for all $x,y$ in $\mathcal{P}$,
$(i)$~$x\leq \overline{x}$, $(ii)$~$x\leq y\Rightarrow \overline{x}\leq \overline{y}$ and
$(iii)$~$\overline{\overline{x}}=\overline{x}$. A general fact of the theory of Heyting algebras is that
performing twice the pseudocomplement operation gives a closure operator. Thus, in the specific case of Dyck algebras,
given a path $P$, its closure $\overline{P}=\sim \sim \! P$ is obtained by
turning each of its factors into the unique pyramid greater than it and having the same number of steps.

\bigskip

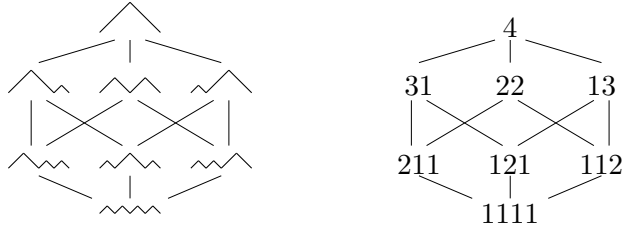
\begin{figure}[h!]
\begin{center}
\begin{tikzpicture}
\begin{scope}[scale=0.1]

\draw (12,0) -- (13,1);
\draw (13,1) -- (14,0);
\draw (14,0) -- (15,1);
\draw (15,1) -- (16,0);
\draw (16,0) -- (17,1);
\draw (17,1) -- (18,0);
\draw (18,0) -- (19,1);
\draw (19,1) -- (20,0);

\draw (0,6) -- (2,8);
\draw (2,8) -- (4,6);
\draw (4,6) -- (5,7);
\draw (5,7) -- (6,6);
\draw (6,6) -- (7,7);
\draw (7,7) -- (8,6);

\draw (12,6) -- (13,7);
\draw (13,7) -- (14,6);
\draw (14,6) -- (16,8);
\draw (16,8) -- (18,6);
\draw (18,6) -- (19,7);
\draw (19,7) -- (20,6);

\draw (24,6) -- (25,7);
\draw (25,7) -- (26,6);
\draw (26,6) -- (27,7);
\draw (27,7) -- (28,6);
\draw (28,6) -- (30,8);
\draw (30,8) -- (32,6);

\draw (0,16)-- (3,19);
\draw (3,19) -- (6,16);
\draw (6,16) -- (7,17);
\draw (7,17) -- (8,16);

\draw (12,16)-- (14,18);
\draw (14,18) -- (16,16);
\draw (16,16) -- (18,18);
\draw (18,18) -- (20,16);

\draw (24,16)-- (25,17);
\draw (25,17) -- (26,16);
\draw (26,16) -- (29,19);
\draw (29,19) -- (32,16);

\draw (12,24) -- (16,28);
\draw (16,28) -- (20,24);

\draw[very thin] (11,2) -- (4,5);
\draw[very thin] (16,2) -- (16,5);
\draw[very thin] (21,2) -- (28,5);

\draw[very thin] (3,9) -- (3,15);
\draw[very thin] (5,9) -- (15,15);
\draw[very thin] (15,9) -- (5,15);
\draw[very thin] (17,9) -- (27,15);
\draw[very thin] (27,9) -- (17,15);
\draw[very thin] (29,9) -- (29,15);

\draw[very thin] (4,20) -- (14,23);
\draw[very thin] (16,20) -- (16,23);
\draw[very thin] (28,20) -- (18,23);

\end{scope}

\begin{scope}[scale=0.1]

\draw (66,0) node {$\small{1111}$};
\draw (54,6.5) node {$\small{211}$};
\draw (66,6.5) node {$\small{121}$};
\draw (78,6.5) node {$\small{112}$};
\draw (54,17) node {$\small{31}$};
\draw (66,17) node {$\small{22}$};
\draw (78,17) node {$\small{13}$};
\draw (66,24.6) node {$\small{4}$};

\draw[very thin] (61,2) -- (54,5);
\draw[very thin] (66,2) -- (66,5);
\draw[very thin] (71,2) -- (78,5);

\draw[very thin] (53,9) -- (53,15);
\draw[very thin] (55,9) -- (65,15);
\draw[very thin] (65,9) -- (55,15);
\draw[very thin] (67,9) -- (77,15);
\draw[very thin] (77,9) -- (67,15);
\draw[very thin] (79,9) -- (79,15);

\draw[very thin] (54,20) -- (64,23);
\draw[very thin] (66,20) -- (66,23);
\draw[very thin] (78,20) -- (68,23);

\end{scope}
\end{tikzpicture}
\end{center}
\caption{The Boolean algebra of regular elements of $\mathfrak{D}_4$
and its isomorphic representation in terms of compositions of 4.\label{regular_compositions}}
\end{figure}

We close by noticing that the Boolean algebra structure of regular elements of $\mathfrak{D}_n$
can be naturally described in terms of \emph{compositions}. Indeed, the map which associates
a concatenation of pyramids in $\mathfrak{D}_n$ with the integer composition (of $n$)
whose parts are the heights of the pyramids (read from left to right) is clearly a bijection.
The partial order induced by $\mathfrak{D}_n$ on the subset of its regular elements can be
translated along such a bijection into the so called \emph{refinement order} on compositions
of $n$, whose covering relation is defined as follows:
a composition $\lambda$ is covered by a composition $\eta$ when $\eta$ is obtained from $\lambda$
by summing two consecutive parts (see Figure \ref{regular_compositions}). These Boolean algebras on compositions have occasionally
surfaced in the literature, see for instance \cite{AS,BLvW,EJ}.

\section{The logic of subintervals}\label{logic}

The aim of this section is to give a logic-theoretic interpretation of Dyck algebras.
More specifically, it turns out that
Dyck algebras provide the natural algebraic counterpart of a special sort of intuitionistic logics,
which are more precisely a certain class of interval temporal logics.

\bigskip

Let $\mathcal{T}_n =\{ t_1 ,t_2 ,\ldots ,t_n \}$ be a finite linearly
ordered set, with $t_1 <t_2 <\cdots <t_n$. The elements of
$\mathcal{T}_n $ will be sometimes called \emph{time states}.
Denote by $Int(\mathcal{T}_n )$ the set of all \emph{intervals}
of $\mathcal{T}_n$, i.e. $I\in Int(\mathcal{T}_n )$ when there
exist $t_i ,t_j \in \mathcal{T}_n $ such that $I=[t_i ,t_j ]=\{
t\; |\, t_i \leq t\leq t_j \}$. In the following we will consider
$Int(\mathcal{T}_n )$ partially ordered by inclusion.

\bigskip

Next we define a set of propositions in a recursive fashion, as usual.
We point out that the logic we are going to describe
is related to the Halpern-Shoham logic \cite{HS},
which is one of the logics of time intervals.
In particular, the propositional logic of interest to us appears to be
intimately related to the fragment of the Halpern-Shoham logic in which
a single modal operator is considered, namely the so-called operator
``during". A paper dealing with this fragment is \cite{MPS}, where
the authors show that it is decidable over finite linear orders.
We also remark that, on the other hand, in \cite{MM} a strictly related fragment
is shown to be undecidable over discrete structures.

The set of propositions $ITL_n$ is defined as follows, by means of the usual connectives:


\begin{itemize}

\item $\bot ,\top \in ITL_n$; for all $1\leq i\leq
n$, $\varepsilon_i \in ITL_n$ (the $\varepsilon_i$'s
are the \emph{propositional variables});

\item if $\varphi ,\psi \in ITL_n$, then $\varphi
\vee \psi , \varphi \wedge \psi , \varphi \rightarrow \psi , \neg
\varphi \in ITL_n$.

\end{itemize}

We give an \emph{interval-based semantics}, for which each
proposition $\varphi$ can be true or false depending on how it is
evaluated on a specific interval $I\in Int(\mathcal{T}_n )$. More
formally, if we denote by $\mathbf{2}^A$ the set of all maps from a
set $A$ to the set $\mathbf{2}=\{ 0,1\}$, we define a map $v$ as
follows:

\begin{eqnarray*}
v&:&ITL_n\longrightarrow \mathbf{2}^{Int(\mathcal{T}_n )}\\
&:&\varphi \longmapsto v_{\varphi}:Int(\mathcal{T}_n
)\longrightarrow \{ 0,1\}
\end{eqnarray*}
where $v_{\varphi}(I)=0$ (resp., 1) if $\varphi$ is false (resp.,
true) when evaluated on the interval $I$. In the following we will
usually write $\varphi (I)$ in place of $v_{\varphi}(I)$. In
particular, we say that $\varphi$ is \emph{valid} when $\varphi
(I)=1$ for all $I\in Int(\mathcal{T}_n )$.

\bigskip

Thus we have a general \emph{evaluation map} $v$, which associates
with every proposition $\varphi$ a specific \emph{valuation}
$v_{\varphi}$ which says on which intervals $\varphi$ is true. The
behavior of valuations with respect to connectives is defined as
usual. More precisely:

\begin{itemize}

\item $(\varphi \vee \psi )(I)=1$ whenever $\varphi (I)=1$ or
$\psi (I)=1$;

\item $(\varphi \wedge \psi )(I)=1$ whenever $\varphi (I)=1$ and
$\psi (I)=1$;

\item $(\neg \varphi )(I)=1$ whenever $\varphi (I)=0$;

\item $(\varphi \rightarrow \psi )(I)=1$ whenever holds: if $\varphi
(I)=1$ then $\psi (I)=1$.

\end{itemize}

Moreover, concerning propositional variables, we define
$\varepsilon_i (I)$ to be true if and only if $I=[t_i ,t_i ]=\{ t_i \}$. We
have therefore all that we need to evaluate any proposition
$\varphi \in ITL_n$.

\bigskip

Notice that, at this point, the partial order structure of
$Int(\mathcal{T}_n )$ does not play any role. We now introduce two
new connectives whose semantics instead depend on such partial
order. These connectives are denoted by $\Box$ and $\lozenge$, and
their semantics is defined as follows:

\begin{itemize}

\item $(\Box \varphi )(I)=1$ when, for all intervals
$J\subseteq I$, $\varphi (J)=1$;

\item $(\lozenge \varphi )(I)=1$ when there exists an
interval $J\subseteq I$ such that $\varphi (J)=1$.

\end{itemize}

Notice that $\Box$ is ``idempotent", in the sense that,
for all intervals $I$, $(\Box \Box \varphi)(I)=(\Box \varphi )(I)$.

\bigskip

We are now ready to describe the subset of $ITL_n$
which will be relevant to us. Define $\Theta_n =\{ \varphi \in
ITL_n\; |\; \varphi \rightarrow \Box \varphi
\textnormal{ is valid} \}$. Intuitively, this means that, if
$\varphi$ is true in $I$, then $\varphi$ is true in all
subintervals of $I$.

\bigskip

We remark here that, from a purely logic-theoretic point of view,
the construction of the set $\Theta_n$ can be suitably described
in the framework of modal companions of an superintuitionistic logic,
see for instance \cite{CZ}. However, the main goal of this section
is to provide a combinatorial description of the logic of $\Theta_n$,
which we believe to be new.

\bigskip

As a subset of $ITL_n$, it is not clear a priori if
$\Theta_n$ is interesting from a semantic point of view. We will
now clarify this point, by showing that $\Theta_n$ is closed with
respect to some, but not all, of the classical connectives.

\begin{prop} If $\varphi ,\psi \in \Theta_n$, then $\varphi \wedge
\psi, \varphi \vee \psi \in \Theta_n$.
\end{prop}

\emph{Proof.}\quad Given $I\in Int(\mathcal{T}_n )$, suppose that
$(\varphi \wedge \psi )(I)=1$, that is $\varphi (I)=\psi (I)=1$.
Since $\varphi ,\psi \in \Theta_n$, we have that, for all
intervals $J\subseteq I$, it is $\varphi (J)=\psi (J)=1$, which
means that $(\Box (\varphi \wedge \psi ))(I)=1$, i.e. $\varphi
\wedge \psi \in \Theta_n$.

Similarly, if we suppose that $(\varphi \vee \psi )(I)=1$, we then
have that $\varphi (I)=1$ or $\psi (I)=1$. Assume, for instance,
that $\varphi (I)=1$. Then, for all intervals $J\subseteq I$, it
is $\varphi (J)=1$, which implies $(\varphi \vee \psi)(J)=1$. We
can thus conclude that $(\Box (\varphi \vee \psi ))(I)=1$, i.e.
$\varphi \vee \psi \in \Theta_n$.\cvd

\begin{prop} $\Theta_n$ is not closed with respect to $\neg$, that is
there exists a proposition $\varphi \in \Theta_n$ such that $\neg
\varphi \notin \Theta_n$.
\end{prop}

\emph{Proof.}\quad Consider the proposition $\varphi
=\varepsilon_1 \vee \varepsilon_2$, and take the interval $I=\{
t_1 ,t_2 \}$. We have clearly $\varphi (I)=0$, and so $(\neg
\varphi )(I)=1$. Now let $J=\{ t_1 \} \subseteq I$: we then get
$\varphi (J)=1$. Therefore we have found an interval $J\subseteq
I$ such that $(\neg \varphi )(J)=0$, which implies that
$(\Box(\neg \varphi ))(I)=0$. We can thus conclude that $((\neg
\varphi )\rightarrow \Box (\neg \varphi ))(I)=0$, as desired.
Notice that this argument clearly works for any proposition of the
type $\varepsilon_i \vee \varepsilon_{i+1}$.\cvd

\begin{prop} $\Theta_n$ is not closed with respect to $\rightarrow$, that is
there exist propositions $\varphi ,\psi \in \Theta_n$ such that
$\varphi \rightarrow \psi \notin \Theta_n$.
\end{prop}

\emph{Proof.}\quad This proposition can be seen as a corollary of
the previous one, since it is not difficult to prove that, for any
interval $I$, $(\neg \varphi )(I)=(\varphi \rightarrow \bot )(I)$.
However we will explicitly provide an example not of that form.

Given $\varphi =\varepsilon_1 \vee \varepsilon_2$ and $\psi
=\varepsilon_2$, we clearly have that $\varphi ,\psi \in
\Theta_n$. Now, given $I=\{ t_2 ,t_3 \}$, we have $\varphi (I)=0$,
hence $(\varphi \rightarrow \psi )(I)=1$. Set $J=\{ t_2 \}
\subseteq I$, we get $\varphi (J)=1$ and $\psi (J)=0$, that is
$(\varphi \rightarrow \psi )(J)=0$. What we have proved so far is
that there is an interval $I$ such that $(\varphi \rightarrow \psi
)(I)=1$ having a subinterval $J$ for which $(\varphi \rightarrow
\psi )(J)=0$. The very last statement (the one concerning $J$)
means that $(\Box (\varphi \rightarrow \psi ))(J)=0$. Therefore we
can conclude that $((\varphi \rightarrow \psi )\rightarrow (\Box
(\varphi \rightarrow \psi )))(I)=0$, and so $\varphi \rightarrow
\psi \notin \Theta_n$.\cvd

The facts that we have recorded so far tell us that the
connectives $\vee$ and $\wedge$ have a nice behavior inside
$\Theta_n$; the same cannot be said for the connectives $\neg$ and
$\rightarrow$. We now define two new
connectives $\sim$ and $\rightsquigarrow$ which can afford better
notions of negation and implication inside $\Theta_n$.

\bigskip

Given an interval $I$ of $\mathcal{T}_n $, we define the semantics
of $\sim$ and $\rightsquigarrow$ as follows:

\begin{itemize}

\item $(\sim \! \varphi )(I)=1$ whenever $\forall J\subseteq I$,
$\varphi (J)=0$;

\item $(\varphi \rightsquigarrow \psi)(I)=1$ whenever $\forall J\subseteq I$,
if $\varphi (J)=1$, then $\psi (J)=1$.


\end{itemize}

Thus, roughly speaking, we say that $\sim \! \varphi$ is true on $I$
whenever $\varphi$ is false on all subintervals of $I$, and that
$\varphi \rightsquigarrow \psi$ is true on $I$ whenever $\psi$ is true
on all subintervals of $I$ on which $\varphi$ is true.
We will call $\sim$ and $\rightsquigarrow$
\emph{pseudonegation} and \emph{pseudoimplication}, respectively.

\bigskip

Observe that the semantics of pseudonegation and pseudoimplication
can be described in terms of classical negation and implication
and the connectives $\Box$ and $\lozenge$. In fact, for any interval $I$,
$(\sim \! \varphi )(I)=(\neg \lozenge \varphi)(I)=(\Box \neg \varphi )(I)$
and $(\varphi \rightsquigarrow \psi)(I)=(\Box (\varphi \rightarrow \psi))(I)$.
Moreover, as an immediate consequence of the definitions, we have
$(\sim \! \varphi)(I)=(\varphi \rightsquigarrow \bot )(I)$.


\bigskip
It is an easy task (and so we leave it to the reader) to prove that,
if $\varphi ,\psi \in \Theta_n$, then
$\sim \! \varphi$, $\varphi ~ \rightsquigarrow ~ \psi \in \Theta_n$.
We now show that pseudonegation has the typical behavior
of an intuitionistic negation.

\begin{prop} Given $\varphi \in \Theta_n$ and $I\in Int(\mathcal{T}_n )$,
if $\varphi (I)=1$, then $(\sim \sim \! \varphi )(I)=1$.
The converse, however, does not hold in general.
\end{prop}

\emph{Proof.}\quad We observe that $(\sim \sim \! \varphi )(I)=1$ if
and only if, for all intervals $J\subseteq I$, there exists an
interval $K\subseteq J$ such that $\varphi (K)=1$. Since $\varphi
\in \Theta_n$, if we suppose that $\varphi (I)=1$, then we have
that, for all intervals $J\subseteq I$, $\varphi (J)=1$, hence the
thesis follows.

To show that the converse does not hold in general, consider the
proposition $\varphi =\varepsilon_1 \vee \varepsilon_2$ and the
interval $I=\{ t_1 ,t_2 \}$. We immediately see that $\varphi
(I)=0$. Moreover, the fact that $(\sim \sim \! \varphi )(I)=1$ is
equivalent to the fact that, for all intervals $J\subseteq \{t_1
,t_2 \}$, there exists an interval $K\subseteq J$ such that
$(\varepsilon_1 \vee \varepsilon_2 )(K)=1$. It is now easy to
realize that the last statement is true.\cvd

\begin{prop} Given $\varphi \in \Theta_n$ and $I\in
Int(\mathcal{T}_n )$, $(\sim \! \varphi )(I)=1$ if and only if $(\sim \sim \sim \! \varphi )(I)=1$.
\end{prop}

\emph{Proof.}
\begin{itemize}

\item[$\Rightarrow$)] 
This is a special case of the previous proposition.

\item[$\Leftarrow$)] Suppose that $(\sim \sim \sim \! \varphi
)(I)=1$, then we have that, for all intervals $J\subseteq I$,
$(\sim \sim \! \varphi )(J)=0$. Thanks to the previous proposition,
this implies that, for all intervals $J\subseteq I$, $\varphi
(J)=0$, that is $(\sim \! \varphi )(I)=1$, as required.\cvd

\end{itemize}

We are now ready to show that pseudonegation and pseudoimplication
are the ``right connectives" in order to describe the Heyting
algebra structure of Dyck paths. Given $\varphi ,\psi \in
\Theta_n$, we say that $\varphi$ and $\psi$ are \emph{equivalent}
when $v(\varphi )=v(\psi )$. In this case we write $\varphi
\vDashv \psi$. It is now left to the reader to show that $\vDashv$
is an equivalence relation on $\Theta_n$ which preserves $\vee
,\wedge ,\rightsquigarrow ,\sim$; this means that, denoting with
$\star$ any of the above mentioned binary connectives, if
$\varphi_1 ,\varphi_2 ,\psi_1 ,\psi_2 \in \Theta_n$ are such that
$\varphi_1 \vDashv \varphi_2$ and $\psi_1 \vDashv \psi_2$, then
$\varphi_1 \star \psi_1 \vDashv \varphi_2 \star \psi_2$ (and a
similar fact holds for the unary connective $\sim$). Thus we can
endow $\Theta_n /{\vDashv}$ with the distributive lattice structure
in which $\vee$ and $\wedge$ are well-defined on equivalence classes
thanks to the above considerations. Denote with $[\Theta_n ]$ the
resulting distributive lattice. Thus, for instance, given
$\varphi ,\psi \in \Theta_n$, denoting with $[\varphi ],[\psi ]\in
\Theta_n /{\vDashv}$ the associated equivalence classes,
in $[\Theta_n ]$ we have that $[\varphi ]\vee [\psi ]=[\varphi \vee \psi ]$,
$[\varphi ]\wedge [\psi ]=[\varphi \wedge \psi ]$ and
$[\varphi ]\rightsquigarrow [\psi ]=[\varphi \rightsquigarrow \psi ]$.
Our next goal is to show that the canonical Heyting algebra structure
on $[\Theta_n ]$ is given precisely by the pseudoimplication operation
$\rightsquigarrow$.

\begin{prop}\label{canonicalheyting} For every $\varphi ,\psi \in \Theta_n$, we have:
$$
[\varphi ]\rightsquigarrow [\psi ] =\bigvee \{ [\alpha ]\in \Theta_n /{\vDashv}
\; |\; [\varphi ] \wedge [\alpha ]\leq [\psi ] \} .
$$

In other words, $\rightsquigarrow$ is the relative pseudocomplement operation
in the canonical Heyting algebra structure of $[\Theta_n ]$.
\end{prop}

\emph{Proof.}\quad We start by observing that the partial order relation $\leq$
associated with the lattice structure of $[\Theta_n ]$ can be described as follows:
$[\varphi ]\leq [\psi ]$ whenever $\varphi (I)\leq \psi (I)$, for all intervals $I$
(which means that, if $\varphi (I)=1$, then $\psi (I)=1$; this is the usual partial order
derived from an algebra of propositions). The reader is invited to see that
$\leq$ is well defined since, if the above condition is satisfied,
then the same inequalities hold when $\varphi$ and $\psi$ are replaced by $\varphi '$ and $\psi '$,
for any $\varphi '\in [\varphi ]$, $\psi '\in [\psi ]$.

Now suppose that $S=\{ [\alpha ]\in \Theta_n /{\vDashv}
\; |\; [\varphi ] \wedge [\alpha ]\leq [\psi ] \} =\{ [\alpha_1 ],[\alpha_2 ],\ldots ,
,[\alpha_r ]\}$. Thus we wish to show that $[\varphi \rightsquigarrow \psi ]=
[\alpha_1 \vee \alpha_2 \vee \cdots \vee \alpha_r ]$. The first step will be to prove
that $[\varphi \rightsquigarrow \psi ]\in S$. Indeed, recall that
the propositions $\alpha_i$ are characterized by the fact that
$[\varphi \wedge \alpha_i ]\leq [\psi ]$. Now, given $I\in Int(\mathcal{T}_n )$,
suppose that $(\varphi \wedge (\varphi \rightsquigarrow \psi ))(I)=1$. This implies that
$\varphi (I)=1$. Then, in order to have $(\varphi \rightsquigarrow \psi )(I)=1$,
necessarily $\psi (I)=1$. This is enough to conclude that
$[\varphi \wedge (\varphi \rightsquigarrow \psi )]\leq [\psi ]$, and so that
$[\varphi \rightsquigarrow \psi ]\in S$, as desired.

To conclude the proof, we will now show that $[\varphi \rightsquigarrow \psi ]$
is an upper bound of $S$, i.e. $[\varphi \rightsquigarrow \psi ]\geq [\alpha_i ]$,
for all $i\leq r$. To this aim, suppose that $\alpha_i (I)=1$, for some interval $I$;
it will be enough to show that $(\varphi \rightsquigarrow \psi)(I)=1$.
Given $J\subseteq I$ such that $\varphi (J)=1$, then we also have $\alpha_i (J)=1$
(since $\alpha_i \in \Theta_n$), and so $(\varphi \wedge \alpha_i )(J)=1$, hence
$\psi (J)=1$. We have thus shown that $(\varphi \rightsquigarrow \psi)(I)=1$,
as desired.\cvd

As usual, to avoid heavy notations, the whole Heyting algebra structure on the set
$[\Theta_n ]$ will simply be denoted $[\Theta_n ]$.
The next lemma is crucial in the proof of our main theorem.

\begin{lemma}\label{CNF} For any $\varphi \in ITL_n$, set $\overline{\varphi}=\sim \sim \! \varphi$.
Given an interval $I$ of $[n]$, set $\varepsilon_I =\overline{\bigvee_{i\in I}\varepsilon_i }$.
Then, for any $\varphi \in \Theta_n$, there exists an antichain of intervals
$I_1 ,I_2 ,\ldots ,I_r$ of $[n]$ such that
$$
\varphi \vDashv \varepsilon_{I_1}\vee \varepsilon_{I_2}\vee
\cdots \vee \varepsilon_{I_r}.
$$
Moreover, when the intervals are listed in increasing order of their minima,
the above one is the unique proposition of that form equivalent to $\varphi$.
\end{lemma}

\emph{Proof.}\quad Fix $\varphi \in \Theta_n$. Denote with $\mathcal{I}\subseteq Int(\mathcal{T}_n )$
the set of all maximal intervals such that $\varphi (I)=1$ (where ``maximal" is intended
with respect to the inclusion order). By construction,
any two elements 
of $\mathcal{I}$
are incomparable; in particular, no two intervals in $\mathcal{I}$ can have either of the two endpoints in common.
Totally order the elements of $\mathcal{I}=\{ I_1 ,I_2 ,\ldots I_r \}$ with respect to their smallest elements
(notice that we would obtain the same total order if we do the same with respect to the greatest elements).
Moreover, identify each element $t_i \in \mathcal{T}_n$ with its index $i\in [n]$. In this way,
we have that $\mathcal{I}\subseteq Int([n])$ and, for each $\alpha \leq r$, $I_{\alpha}\in Int([n])$.
Our aim is now to prove
\begin{equation}\label{equivalence}
\varphi \vDashv \bigvee_{1\leq \alpha \leq r}\varepsilon_{I_{\alpha}}.
\end{equation}
Before starting to prove this equivalence, it is convenient to observe the following two facts:
\begin{itemize}
\item for all intervals $J\subseteq I$, $\varepsilon_I (J)=1$;
\item for all intervals $J\nsubseteq I$, $\varepsilon_I (J)=0$.
\end{itemize}
Indeed, given an interval $J\subseteq I$, we have $\varepsilon_I (J)=\sim \sim \! (\bigvee_{i\in I}\varepsilon_i )(J)=1$
if and only if, for all intervals $K\subseteq J$, there exists an interval $M\subseteq K$ such that
$$
\bigvee_{i\in I}\varepsilon_i (M)=1.
$$

The last statement is in fact true: for a given $K\subseteq J$, it is enough to choose an element $\tau \in K$
in order to have $\bigvee_{i\in I}\varepsilon_i (\{ \tau \} )\geq \varepsilon_{\tau }(\{ \tau \} )=1$.

On the other hand, given an interval $J\nsubseteq I$, we have
$\varepsilon_I (J)=\sim \sim \! (\bigvee_{i\in I}\varepsilon_i )(J)=0$ if and only if
there exists an interval $K\subseteq J$ such that, for all intervals $M\subseteq K$,
$$
(\bigvee_{i\in I}\varepsilon_i )(M)=0.
$$

Once again, it is not difficult to see that the last equality is true:
choosing, for instance, $K=J\setminus I$, one immediately realizes that,
for every $i\in I$, $\varepsilon_i (M)=0$ (since $i\notin M$, and so $M\neq \{ i\}$).

We are now ready to proceed with the announced proof of (\ref{equivalence}).
Given an interval $I$, since the only possible truth values are 0 and 1,
it will be enough to prove what follows:
\begin{itemize}
\item[(i)] if $\varphi (I)=1$, then $\bigvee_{\alpha}\varepsilon_{I_{\alpha}}(I)=1$;
\item[(ii)] if $\varphi (I)=0$, then $\bigvee_{\alpha}\varepsilon_{I_{\alpha}}(I)=0$.
\end{itemize}

Let us prove the two above statements separately.

\begin{itemize}
\item[(i)] Suppose that $\varphi (I)=1$. Then there exists $s$ such that $I_s\in \mathcal{I}$ and $I\subseteq I_s$.
Therefore $\bigvee_{\alpha}\varepsilon_{I_{\alpha}}(I)\geq \varepsilon_{I_s}(I)=1$.
\item[(ii)] Suppose that $\varphi (I)=0$. This means that $I\nsubseteq I_s$, for all $s\leq r$.
Therefore $\varepsilon_{I_s}(I)=0$, for all $s$, hence $\bigvee_{\alpha}\varepsilon_{I_{\alpha}}(I)=0$.\cvd
\end{itemize}

For any given $\varphi \in \Theta_n$, the above lemma provides a canonical form for $\varphi$,
which will be called its \emph{closed disjunctive form} (briefly, \emph{CDF}).

The next theorem is the main result of the present paper.

\begin{teor}\label{representation} The Heyting algebra $\mathfrak{D}_n$ of Dyck paths of
semilength $n$ is isomorphic to the Heyting algebra $[\Theta_{n-1}]$.
\end{teor}

\emph{Proof.}\quad By the previous lemma, we can (and in fact will) identify
each equivalence class of $[\Theta_n ]$ with the unique proposition in CDF
contained in the class. Moreover, we recall that, in $\mathfrak{D}_n$,
the atoms are those paths all of whose factors are hills except for a single pyramid
having exactly 4 steps. If $P$ is an atom of $\mathfrak{D}_n$, we denote with $x_P$
the abscissa of the unique nontrivial peak of $P$, and we call $x_P /2$ the \emph{order} of the atom $P$.

Define the function $f:[\Theta_{n-1}]\rightarrow \mathfrak{D}_n$ as follows:
given pairwise incomparable intervals $I_1 ,I_2 ,\ldots ,I_r \subseteq [n-1]$, set
$f(\varepsilon_{I_1}\vee \varepsilon_{I_2}\vee \cdots \vee \varepsilon_{I_r})$ equal to
the Dyck path $P$ of semilength $n$ whose decomposition into join-irreducibles
$P=P_1 \vee P_2 \vee \cdots \vee P_r$ has cardinality $r$ and is such that, for every $j\leq r$,
the interval of atoms below $P_j$ is made by the atoms of order $i$, for all $i\in I_j$.
We claim that $f$ is a Heyting algebra isomorphism.

We start by showing that $f$ is onto. Indeed, given any Dyck path $P$ in $\mathfrak{D}_n$,
its decomposition into join-irreducibles uniquely determines an antichain of intervals of $[n-1]$,
which is given by the intervals $I_1 ,\ldots I_r$ of the orders of the atoms lying below each join-irreducible.
By construction, the proposition (in CDF) $\varepsilon_{I_1}\vee \cdots \vee \varepsilon_{I_r}$
is mapped by $f$ onto $P$.

Next we prove that $f$ is order-preserving. To this aim, we first give an alternative description of
the partial order of the Heyting algebra $[\Theta_{n-1}]$, based on the CDF representatives of equivalence classes.
Given $\varphi ,\psi$ in $[\Theta_{n-1}]$, suppose that $\varphi =\varepsilon_{I_1}\vee \cdots \vee \varepsilon_{I_r}$
and $\psi =\varepsilon_{J_1}\vee \cdots \vee \varepsilon_{J_s}$, for suitable antichains of intervals in $[n-1]$.
Recall that $\varphi \leq \psi$ if and only if, for all $I\subseteq [n-1]$, $\varphi (I)\leq \psi (I)$.
Our assumptions on $\varphi$ and $\psi$ implies that
$\varphi (I)=1$ if and only if $I\subseteq I_h$, for some $h\leq r$ (and analogously for $\psi$).
Thus we get that $\varphi \leq \psi$ if and only if, for every $h\leq r$, there exists $k\leq s$
such that $I_h \subseteq J_k$. Now suppose that $\varphi \leq \psi$. If $r=s=1$,
then $f(\varphi )=P$ and $f(\psi )=Q$ are join-irreducibles in $\mathfrak{D}_n$, i.e.
they consist of a series of hills and a unique pyramid having at least 4 steps. Saying that
$\varphi \leq \psi$ means in this case that $I_1 \subseteq J_1$, hence the interval of atoms dominated by $P$
is contained in the interval of atoms dominated by $Q$, that is $P\leq Q$.
In the general case, set $f(\varphi )=P=P_1 \vee \cdots \vee P_r$ and
$f(\psi )=Q=Q_1 \vee \cdots \vee Q_s$; if $\varphi \leq \psi$, then, for every $h\leq r$,
there exists $k\leq s$ such that $I_h \subseteq J_k$, which implies that $P_h \leq Q_k$.
From here it follows that $P\leq Q$.

All the above arguments can be reversed, thus showing that $f$ is also order-reflecting,
i.e. that $f(\varphi )\leq f(\psi )$ implies that $\varphi \leq \psi$.

Therefore we have shown that $f$ is onto, order-preserving and order-reflecting. It is known that
this is enough to conclude that $f$ is an order isomorphism. As a consequence, $f$ is also a
lattice isomorphism. Finally, thanks to Proposition \ref{canonicalheyting}, if we consider the
canonical Heyting algebra structure induced by the finite distributive lattice structure, we have that
$f$ is a Heyting algebra isomorphism between $\mathfrak{D}_n$ and $[\Theta_{n-1}]$, as desired.\cvd

\section{Posets of intervals}\label{poset}


The results of the previous sections suggest that every element of a Dyck algebra
can be described by means of the underlying Heyting algebra structure
together with a natural linear order structure on the set of the atoms of the algebra.
Below we will try to clarify this statement.

\bigskip

Given a Dyck path $P$, denote with $\overline{P}$ its Heyting algebra closure, that is
$\overline{P}=\sim \sim \! P$. The set of atoms of a Dyck algebra can be given a total order structure
(which has nothing to do with the partial order of the algebra) by declaring an atom $P$
strictly less than another atom $Q$ whenever $x_P <x_Q$
(we refer to the notation introduced in the proof of theorem \ref{representation} for the order of an atom). 
In this case we will write $P\ll Q$, to avoid confusion with the partial order on Dyck paths.
The (finite) set of atoms of $\mathfrak{D}_n$ will then be denoted $\{ \pi_1 ,\pi_2 ,\ldots ,\pi_{n-1}\}$,
where $\pi_i$ is the atom of order $i$.
As we have already noticed, a join-irreducible path is uniquely determined by the set of atoms lying below it.
Such a set of atoms is obviously an interval with respect to $\ll$.
More specifically, if $P$ is a join-irreducible and $\pi_i ,\pi_{i+1},\ldots ,\pi_{i+j}$
are the atoms below $P$, then $P=\overline{\pi_i \vee \pi_{i+1}\vee \cdots \vee \pi_{i+j}}$.
Summing up, every Dyck path can be expressed (via Birkhoff representation theorem) as
the join of the closure of the join of $\ll$-intervals of atoms.

A further step towards abstraction consists of identifying an interval of atoms of $\mathfrak{D}_n$
with the interval of the orders of such atoms (which is a subset of $[n-1]
$). Thus a Dyck path of semilength $n$
can be identified with a family of incomparable intervals (i.e., an antichain of intervals) of $[n-1]$.
This observation leads to a possible generalization of the approach
we have developed so far for Dyck algebras, which we attempt to sketch in the remainder of this section.

\bigskip

Let $\mathcal{P}$ be a poset and denote with $Int(\mathcal{P})$ the poset of bounded intervals of $\mathcal{P}$
ordered by inclusion. The generic element of $Int(\mathcal{P})$ is then
$[x,y]=\{ z\in \mathcal{P}\; |\; x\leq z\leq y\}$.  We are interested in the set $\mathcal{O}(Int(\mathcal{P}))$
of all down-sets of $Int(\mathcal{P})$. When ordered by inclusion, $\mathcal{O}(Int(\mathcal{P}))$ is a
complete distributive lattice. This kind of lattices is often relevant from a theoretical point of view.
For instance, we recall here that, when $\mathcal{P}$ is locally finite (i.e. every interval of $\mathcal{P}$ is finite),
$\mathcal{O}(Int(\mathcal{P}))$ is isomorphic to the lattice of two-sided ideals of the incidence algebra of $\mathcal{P}$.
This is a crucial fact in showing that two locally finite posets are order-isomorphic if and only if
their incidence algebras are isomorphic (see, for instance, \cite{DRS}).

\begin{lemma} The lattice $\mathcal{O}(Int(\mathcal{P}))$ is atomic (i.e. every element of $\mathcal{O}(Int(\mathcal{P}))$
contains at least one atom), and the set of its atoms is in bijection with $P$.
\end{lemma}

\emph{Proof.}\quad For any $x\in \mathcal{P}$, the interval $[x,x]$ is a minimal element of $Int(\mathcal{P})$
(and every minimal element is of this form). Therefore the set
$\mathcal{A}=\{ \{ \emptyset ,[x,x]\} \subseteq Int(\mathcal{P})\; |\; x\in \mathcal{P}\}$ is the set of atoms of
$\mathcal{O}(Int(\mathcal{P}))$. Since every nonempty down-set of $Int(\mathcal{P})$ contains at least one interval
$I$, if $x\in I$, then obviously $\{ \emptyset ,[x,x]\}$ is contained in the given down-set, which is enough
to conclude.\cvd


The above lemma asserts that there is a natural partial order on the set of atoms of $\mathcal{O}(Int(\mathcal{P}))$
(inherited from the partial order of $\mathcal{P}$), which has of course nothing to do with the inclusion order on $\mathcal{O}(Int(\mathcal{P}))$.
It would be very interesting to deduce properties of the complete distributive lattice $\mathcal{O}(Int(\mathcal{P}))$ from properties of $\mathcal{P}$.
Since lattices of down-sets are completely distributive, they are also Heyting algebras (in the same canonical way as finite lattices are),
thus the same project can be developed for the Heyting algebra structure of $\mathcal{O}(Int(\mathcal{P}))$.
To the best of our knowledge, it seems that this approach to the study of posets of intervals has never been considered before.
To justify it, we now briefly mention some remarkable examples.



\bigskip

\emph{Examples.}

\begin{enumerate}

  \item If $\mathcal{P}$ is a discrete poset (i.e., an antichain), then clearly $Int(\mathcal{P})\simeq \mathcal{P}$, hence
  any element of $\mathcal{O}(Int(\mathcal{P}))$ can be seen as a subset of $\mathcal{P}$. This means that
  $\mathcal{O}(Int(\mathcal{P}))$ is a complete and atomic Boolean algebra.

  \item If $\mathcal{P}$ is totally ordered, then, in the finite case,
  $\mathcal{O}(Int(\mathcal{P}))$ is isomorphic to a Dyck lattice of suitable order, see also
  \cite{FM1}. In case $\mathcal{P}$ is infinite, we obtain a natural infinite analog of Dyck lattices
  which still deserves to be studied.

  \item If $\mathcal{P}$ is a finite Boolean algebra, then $Int(\mathcal{P})$ is the sup-semilattice
  of the nonempty faces of a cube of suitable dimension (see \cite{BO1,BO2}).
  However, the distributive lattice $\mathcal{O}(Int(\mathcal{P}))$ has never been studied;
  a better understanding of its structure, as well as of its logic-theoretic properties as a Heyting algebra,
  is surely desirable. Also, we are not aware of what happens for infinite Boolean algebras.

\end{enumerate}


\section{Combinatorial properties of Dyck paths in terms of atoms of Dyck lattices}\label{combinatorics}


In this final section we will give a glimpse of the potential applications of the general approach
outlined in the previous section in the particular case of Dyck algebras. More specifically,
we will focus on combinatorics, and we will express several combinatorial properties of Dyck paths
in terms of  the Heyting algebra structure of Dyck algebras and the natural linear order $\ll$ on their atoms.

We recall once again that every path of the Dyck algebra $\mathfrak{D}_n$ can be identified with
an antichain of intervals of the totally ordered set $[n-1]$ (namely, the family of
pairwise incomparable intervals each of which represents the indices of the atoms
dominated by a join-irreducible in the decomposition of the path). For instance,
the red Dyck path in Figure \ref{implication} corresponds to the antichain of intervals
$\{ [2,4],[4,5],[6,6],[8,8],[9,9]\}$ of the set $[9]$. For any two such antichains
$\{ I_1 ,\ldots ,I_n \}$ and $\{ J_1 ,\ldots ,J_m \}$, it is $\{ I_1 ,\ldots ,I_n \} \leq \{ J_1 ,\ldots ,J_m \}$
in $\mathfrak{D}_n$ whenever, for every $i\leq n$, there exists $j\leq m$ such that
$I_i \subseteq J_j$ (as we already noticed in the proof of Theorem \ref{representation}).
It is also useful to record an explicit expression for join and meet:
\begin{eqnarray*}\label{meetandjoin}
\{ I_1 ,\ldots ,I_n \} \vee \{ J_1 ,\ldots ,J_m \} &=&\{ I_1 ,\ldots ,I_n ,J_1 ,\ldots ,J_m \} ;\nonumber \\
\{ I_1 ,\ldots ,I_n \} \wedge \{ J_1 ,\ldots ,J_m \} &=&\{ I_i \cap J_j \; |\; i\leq n,j\leq m \} ,\nonumber
\end{eqnarray*}
where in both the r.h.s.'s we tacitly assume to discard all intervals that are not maximal
(this is of course needed in order to get an antichain).
We can also give a description of pseudonegation:
if a path $P$ is represented by the antichain of intervals $\{ I_1 ,\ldots ,I_m \}$, then $\sim \! P$
is represented by the (unique) family of maximal intervals constituting a partition of the set
$[n-1]\setminus (I_1 \cup \cdots \cup I_m )$. Referring to the black path in Figure \ref{pseudocomplement},
its pseudonegation is represented by the antichain of intervals $\{ [4,4],[12,14]\}$ of $[15]$.

\bigskip

We now state and prove a series of propositions which express some important combinatorial parameters
on Dyck paths in terms of the above described ``interval" representation of Dyck paths.
For a classical reference on the enumerative combinatorics of Dyck paths, see the survey article
\cite{D}. Before starting we need to introduce a few notations and definitions.

For a given Dyck path $P\in D_n$, we denote with $\mathcal{F}_P$ the antichain of intervals of $[n-1]$
representing that path. If $\mathcal{F}_P=\{ I_1 ,\ldots ,I_m \}$, then the \emph{cardinality} of
$\mathcal{F}_P$ is $|\mathcal{F}_P|=m$, whereas the \emph{weight} of $\mathcal{F}_P$ is
$\| \mathcal{F}_P\| =|I_1 \cup \cdots \cup I_m |$. Moreover, we say that $I\in \mathcal{F}_P$ is
\emph{internal} when $1,n-1\notin I$; the set of internal intervals of $\mathcal{F}_P$ is denoted
with $\mathcal{F}_P ^*$.

\begin{prop} The number of peaks of a Dyck path $P\in D_n$ is given by
$|\mathcal{F}_P |+\| \mathcal{F}_{\sim \! P}\| -|\mathcal{F}_{\sim \! P}^*|$.
\end{prop}

\emph{Proof.}\quad Each peak of $P$ of height $>1$ represents the contribution of a join-irreducible
in the (unique) expansion of $P$ as a join of join-irreducibles. Since join-irreducibles of $P$
correspond to intervals of $\mathcal{F}_P$, the contributions of these peaks is exactly $|\mathcal{F}_P |$.
As far as peaks at height 1 are concerned (i.e., hills), we observe that
a bunch of $s$ consecutive hills of $P$ corresponds to an internal interval of cardinality $s+1$ of $\mathcal{F}_{\sim \! P}$,
except when the bunch of hills is at the beginning or at the end of the path,
in which cases it corresponds to a noninternal interval of cardinality $s$ of $\mathcal{F}_{\sim \! P}$.
This means that the number of hills of $P$ is $\| \mathcal{F}_{\sim \! P}\| -|\mathcal{F}_{\sim \! P}^*|$,
which concludes the proof.\cvd

A byproduct of the above proof is the following.

\begin{cor} The number of hills of a Dyck path $P\in D_n$ is given by
$\| \mathcal{F}_{\sim \! P}\| -|\mathcal{F}_{\sim \! P}^*|$.
\end{cor}

\begin{prop} The sum of the heights of the peaks of a Dyck path $P\in D_n$ is given by
$\| \mathcal{F}_P \| +|\mathcal{F}_P|+\| \mathcal{F}_{\sim \! P}\| -|\mathcal{F}_{\sim \! P}^*|=
n-1+|\mathcal{F}_P|-|\mathcal{F}_{\sim \! P}^*|$.
\end{prop}

\emph{Proof.}\quad Concerning peaks of height $>1$, we observe that the height of each of them is
the cardinality of the interval which correspond to it minus 1. Thus the contribution to the total heights sum
of such peaks is $\| \mathcal{F}_P \|+|\mathcal{F}_P |$. On the other hand, the sum of the heights of the hills
of $P$ equals the number of hills of $P$, so (from the proof of the previous proposition), their contribution
is given by $\| \mathcal{F}_{\sim \! P}\| -|\mathcal{F}_{\sim \! P}^*|$. Summing up the two quantities we have obtained
gives the desired result.\cvd

\begin{prop} The number of returns of a Dyck path $P\in D_n$ is given by
$\| \mathcal{F}_{\sim \! P}\| +1$.
\end{prop}

\emph{Proof.}\quad The total number of returns of $P$ is given by the number of its hills
plus the number of its nontrivial factors. As we have already proved, the number of hills of $P$ is given by
$\| \mathcal{F}_{\sim \! P}\| -|\mathcal{F}_{\sim \! P}^*|$. Moreover
we observe that the number of nontrivial factors of $P$ is ``approximately equal" to the number of
nontrivial factors of $\sim \! P$. They are indeed equal if and only if $P$ either starts or ends with a hill
(but not both); in this case, $P$ has precisely $|\mathcal{F}_{\sim \! P}|$ nontrivial factors,
and so the total number of returns of $P$ is
$\| \mathcal{F}_{\sim \! P}\| -|\mathcal{F}_{\sim \! P}^*|+|\mathcal{F}_{\sim \! P}|=\| \mathcal{F}_{\sim \! P}\|+1$
(since in this case $\sim \! P$ has precisely one nontrivial factor either at the beginning or at the end,
which corresponds to a single noninternal interval). Otherwise,
$P$ has one more (resp., less) nontrivial factor than $\sim \! P$ if and only if $P$ both starts and ends
with a nontrivial factor (resp., with a hill); in this case $P$ has precisely
$|\mathcal{F}_{\sim \! P}|+1$ (resp., $|\mathcal{F}_{\sim \! P}|-1$) nontrivial factors,
and so the total number of returns of $P$ is
$\| \mathcal{F}_{\sim \! P}\| -|\mathcal{F}_{\sim \! P}^*|+|\mathcal{F}_{\sim \! P}|+1$
(resp., $\| \mathcal{F}_{\sim \! P}\| -|\mathcal{F}_{\sim \! P}^*|+|\mathcal{F}_{\sim \! P}|-1$),
which equals $\| \mathcal{F}_{\sim \! P}\|+1$ (the reader is invited to check all the details).\cvd

All the results illustrated so far concern statistics which can be directly expressed in terms of
global parameters. We give below a few simple examples in which it is necessary to take into account
some local information. The last example is especially interesting, being an instance of a kind of
``pattern occurrence" statistic. Since the proofs are quite easy, we leave most of them to the reader.
Recall that the ``interval" representation of a generic Dyck path $P$ is written $\{ I_1 ,\ldots ,I_m \}$,
where each $I_i$ is an interval of $[n-1]$, and the intervals are listed in increasing order of their minima.
Moreover, we say that two consecutive intervals $I_i$ and $I_{i+1}$ are \emph{distanced} when
$\max I_i < \min I_{i+1}-1$.

\begin{prop} The height of the first peak of a Dyck path $P\in D_n$ is given by
\begin{eqnarray*}
\left\{ \begin{array}{ll}
|I_1 |+1 & \textnormal{, if $1\in I_1$} \\
1 & \textnormal{, otherwise}
\end{array}\right. .
\end{eqnarray*}
\end{prop}

\begin{prop} The number of peaks before the first return of a Dyck path $P\in D_n$ is given by
\begin{eqnarray*}
\left\{ \begin{array}{ll}
\max \{ k\; |\; \textnormal{$I_{i-1}$ and $I_i$ are not distanced, for all $i\leq k$}\} & \textnormal{, if $1\in I_1$} \\
1 & \textnormal{, otherwise}
\end{array}\right. .
\end{eqnarray*}
\end{prop}

\begin{prop} The number of occurrences of the (consecutive) factor $duu$ in a Dyck path $P\in D_n$
is given by
$$
|\{ i\leq n-1\; |\; \textnormal{either $I_{i-1}$ and $I_i$ are distanced or $|I_i \setminus I_{i-1}|>1$}\}.
$$
\end{prop}

\emph{Proof.}\quad Each occurrence of $duu$ in $P$ corresponds to the occurrence of a valley not immediately followed by a peak.
For any such valley we have two distinct possibilities.
If the valley is not on the $x$-axis, then it corresponds to a transition between two consecutive join-irreducibles such that
the rightmost one dominates at least two atoms which are not dominated by the leftmost one.
In terms of the ``interval" representation of the path, this corresponds to a consecutive pair of
non-distanced intervals $I_{i-1}$ and $I_i$ such that $|I_i \setminus I_{i-1}|>1$.
On the other hand, if the valley lies on the $x$-axis, then it is immediately followed by a nontrivial factor,
and the first interval $I_i$ corresponding to such a factor is clearly distanced from the previous one $I_{i-1}$.
\cvd

%
%
%
%
%

\section{Further work}\label{open}

As already illustrated in section \ref{poset}, the case of Dyck algebras investigated here
is just an instance of a more general situation. The study of the complete distributive lattices
(Heyting algebras) of the down-sets of the poset of intervals of a generic poset is a totally unexplored subject,
which seems interesting to be pursued both from the algebraic and the logic-theoretic point of view.
We remark that the relevance of posets of intervals in certain logical framework has already been noticed,
see \cite{CM}.

In particular, the case in which the starting poset is a Boolean algebra (example 3 in section \ref{poset})
is related to the logic of the $n$-cube, initiated in \cite{RM} and recently explored in \cite{Mun}.

\bigskip

It would be nice to have a purely algebraic characterization of Dyck lattices and of Dyck algebras.
Even if they are not a variety (in the sense of universal algebra), they show some interesting features.
For instance, Dyck lattices are projective distributive lattices (this follows from a result of Balbes \cite{B},
which asserts that a finite distributive lattice is projective if and only if the poset of its join-irreducibles
is a meet-semilattice).

\bigskip

It is natural to replace Dyck paths with other families of paths. The first, obvious candidates
are Motzkin and Schr\"oder paths. In both cases, the analogous posets are distributive lattices too,
so an investigation of the associated Heyting algebra structures and of their logic-theoretic interpretations
can be done along similar lines.

\bigskip

Is it possible to find analogous results for other fragments of the Halpern-Shoham logic?
More specifically, are there similar combinatorial descriptions when the underlying order of time instants is a finite total order?

\bigskip

{\footnotesize {\bf Acknowledgment.}\quad The author wishes to warmly thank Daniele Mundici,
for many stimulating discussions on the topics presented here, as well as for several remarks and comments
on an earlier draft, which have contributed to significantly improve it.}


%
%
%
%
%
%
%
%
%
%

\end{document}